\documentclass[11pt, a4paper, xy]{article}
\usepackage[all]{xy}
\usepackage{amssymb,amsmath}

\usepackage{amsfonts}
\usepackage{amscd}

\usepackage{url}

\makeatletter
\def\url@leostyle{%
  \@ifundefined{selectfont}{\def\UrlFont{\sf}}{\def\UrlFont{\small\ttfamily}}}
\makeatother
\input amssym.def
\input amssym

\newtheorem{thm}{Theorem}[section]
\newtheorem{lem}[thm]{Lemma}
\newtheorem{prop}[thm]{Proposition}
\newtheorem{cor}[thm]{Corollary}

\newtheorem{clm}[thm]{Claim}

\newtheorem{defn}[thm]{Definition}

\newtheorem{nrmk}[thm]{Remark}
\newtheorem{expl}[thm]{Example}

\newcommand{\qed}{\hfill $\Box$ \vspace{.5cm}}
\newcommand{\pf}{{\bf Proof. }}

\title {The Lefschetz coincidence theorem in o-minimal expansions of fields}

\author {M\'{a}rio J. Edmundo \thanks{With partial support from the FCT 
(Funda\c{c}\~ao para a Ci\^encia e Tecnologia), program POCTI (Portugal/FEDER-EU). 
}\\ 
{\small CMAF Universidade de Lisboa}\\ 
{\small Av. Prof. Gama Pinto 2}\\
{\small 1649-003 Lisboa, Portugal}\\
{\small edmundo@cii.fc.ul.pt}\\
{\small and}\\
{Arthur Woerheide}\\
{\small 6203 S. Evans Ave.}\\
{\small Chicago, IL 60637}\\
{\small U.S.A}\\
\\}
\date{June 29, 2009}

\newcommand{\into}{\longrightarrow}

\renewcommand{\hat}{\widehat}
\renewcommand{\tilde}{\widetilde}
\renewcommand{\bar}{\overline}

\newcommand{\NN}{\mathbb{N}}
\newcommand{\ZZ}{\mathbb{Z}}
\newcommand{\QQ}{\mathbb{Q}}

\newcommand{\N}{\mbox{$\cal N$}}

\begin{document}

\maketitle

\begin{abstract}
In this paper we prove the Lefschetz coincidence theorem in o-minimal expansions of 
fields using the o-minimal singular homology and cohomology.
\end{abstract}

\newpage

\begin{section}{Introduction}\label{section introduction}
We work over an  o-minimal expansion $\N$$=(N,0,1,<,+,\cdot ,\dots)$ of a 
real closed field $N$. Definable means $\N$-definable (possibly with 
parameters). As it is well known, o-minimal structures are a wide ranging 
generalization of semi-algebraic and sub-analytic geometry. Good 
references on o-minimality are, for example, the book \cite{vdd} by 
van den Dries and the notes \cite{c} by Michel Coste. For semi-algebraic geometry relevant to this paper the reader should consult  the work by Delfs and Knebusch  such as \cite{dk} and the book \cite{BCR} by Bochnak, Coste and M-F. Roy.
 
The goal of the paper is to present the proof of the following o-minimal version of the 
Lefschetz coincidence theorem.

\begin{thm}[Lefschetz coincidence theorem]\label{thm lct}
Let $X$ and $Y$ be orientable, definably compact definable manifolds of dimension $n$. Suppose that  $f,g:X\into Y$ are continuous definable maps whose  Lefschetz coincidence 
number  is nonzero. Then there is $x\in X$ such that $f(x)=g(x)$.
\end{thm}

This result implies an o-minimal Lefschetz fixed point theorem for definable continuous maps on orientable, definably compact definable manifolds as in \cite{bo}. For a more general o-minimal Lefschetz-Hopf fixed point theorem generalizing  Brumfiel's Hopf fixed point theorem for semi-algebraic maps in \cite{br} see \cite{e2}.

Our  proof of o-minimal Lefschetz coincidence theorem above follows the proofs of its 
topological analogue (\cite{g}, \cite{ghv}). The only difficulty is the o-minimal relative Poincar\'e duality theorem for orientable definable manifolds (Theorem \ref{thm poincare duality}) which is handled by replacing classical arguments such as compactness (\cite{d}) by the definable triangulation theorem (\cite{vdd}) and the 
existence of coverings by definable proper sub-balls (\cite{bo}, \cite{e3}, \cite{w}). With this result  available in the o-minimal setting, we develop in Subsection 
\ref{subsection thom, lefschetz and  euler classes} the o-minimal analogue of part of the classical theory of  Thom,  Lefschetz and  Euler classes as in (\cite{d}) and  prove in a rather classical and algebraic way the Lefschetz coincidence theorem  in Subsection \ref{subsection lefschetz fixed point theorem}.

\end{section}

\begin{section}{O-minimal (co)homology}
\label{section preli}

For o-minimal expansions of real closed fields, Woerheide (\cite{Wo}) 
constructs o-minimal singular homology $(H_*,d_*)$ with coefficients 
in ${\ZZ}$ satisfying o-minimal Eilenberg-Steenrod homology axioms 
(the analogues of the classical Eilenberg-Steenrod axioms for the 
category of definable sets with definable continuous maps).

The definition of o-minimal singular homology is quite easy, but the 
verification of the  axioms is very difficult as we now explain. Given 
a definable set $X$ we consider, for each $m\geq 0$, the abelian  group 
$S_m(X)$ freely generated by the singular definable simplices 
$\sigma :\Delta ^m\into X$, where 
$\Delta ^m=\{(t_0,\dots ,t_m)\in N^{m+1}:\sum _it_i=1,\,t_i\geq 0\}$ 
is the standard $m$-dimensional simplex. The boundary operator
$\partial _{m+1}:S_{m+1}(X)\into S_m(X)$ (morphism of degree $-1$) is 
defined as in the classical case making $S_*(X)$ a free chain complex. 
Also, a definable continuous map $f:X\into Y$ induces a chain map 
$f_{\sharp }:S_*(X)\into S_*(Y)$ (i.e., a morphism of degree zero 
satisfying $f_{\sharp }\circ \partial _*=\partial _*\circ f_{\sharp }$). 
Similarly one defines the definable singular chain complex of pairs of 
definable sets $A\subseteq X$ by $S_*(X,A)=S_*(X)/S_*(A)$. The graded 
group $H_*(X)$ is defined as the homology of the chain complex 
$S_*(X)$. Similarly $H_*(X,A)$ is the homology of
$S_*(X,A)$. A definable continuous map $f:X\into Y$ induces a homomorphism
$f_*:H_*(X)\into H_*(Y)$ of graded groups (via $f_{\sharp }$). In the 
same way, a definable continuous map $f:(X,A)\into (Y,B)$ (i.e., a 
definable continuous map $f:X\into Y$ such that $f(A)\subseteq B$) 
induces a homomorphism
$f_*:H_*(X,A)\into H_*(Y,B)$ of graded groups. 

\begin{thm}[Woerheide \cite{Wo}]\label{thm woer}
The o-minimal homology functor $H_*$ constructed above satisfies the 
o-minimal Eilenberg-Steenrod axioms:

{\bf Homotopy Axiom}. If $f,g:(X,A)\into (Y,B)$ are definable maps and 
there is a definable homotopy between $f$ and $g$, then \[ f_*=g_*:
H_n(X,A)\into H_n(Y,B) \] for all $n\in \NN$. 

{\bf Exactness Axiom}. For $A\subseteq X$ definable sets if 
$i:(A,\emptyset )\into (X,\emptyset )$ and $j:(X,\emptyset )\into  (X,A)$ 
are the inclusions, then we have a natural exact sequence 
\[ \cdots \into H_n(A,\emptyset)\stackrel{i_*}
{\rightarrow} H_n(X,\emptyset)\stackrel{j_*}{\rightarrow}H_n(X,A) 
\stackrel{d_n}{\rightarrow}H_{n-1}(A,\emptyset)\into \cdots .
\]

{\bf Excision Axiom}. For every pair $A\subseteq X$ of definable sets 
and every definable open  subset $U$ of $X$ such that 
$\bar{U}\subseteq \mathring{A}$, the inclusion $(X-U,A-U)\into (X,A)$ 
induces isomorphisms 
\[ H_n(X-U,A-U) \into H_n(X,A) \] for all $n\in \NN$.

{\bf Dimension Axiom}. If $X$ is a one point set, then 
$H_n(X,\emptyset)=0$ for all $n\neq 0$ and $H_0(X)={\ZZ}$. 
\end{thm}

Woerheide's result is based on a definable triangulation theorem 
(\cite{vdd}) and on the method of acyclic models from homological 
algebra and is rather complicated due to the fact that, 
in arbitrary o-minimal expansions of fields, the classical simplicial 
approximation theorem, the method of repeated barycentric subdivisions 
and the Lebesgue number property for a standard simplex 
$\Delta ^n$ fail.

We make now a few comments comparing the classical proof of the excision axiom and Woerhiede proof of the o-minimal excision axiom. 

For $z\in \tilde{S}_*(X)$ with $z=\sum _{j=1}^la_j\alpha _j$ we have a chain map 
$$z_{\sharp }:\tilde{S}_*(\Delta ^n)\into \tilde{S}_*(X):\beta \mapsto z_{\sharp }\beta =\sum _{i,j}a_jb_i(\alpha _j\circ \beta _i)$$
where $\beta =\sum _{i=1}^kb_i\beta _i$. 

Let $X$ be a definable set. The {\it barycentric subdivision} 
$${\rm Sd}_n:\tilde{S}_n(X)\into \tilde{S}_n(X)$$ 
is defined as follows: for $n\leq -1$, ${\rm Sd}_n$ is the trivial homomorphism, ${\rm Sd}_{-1}$ is the identity and, for $n\geq 0$, we set 
$${\rm Sd}_n(z)=z_{\sharp }(b_n.{\rm Sd}_{n-1}\partial 1_{\Delta ^n})$$ 
where $b_n$ is the barycentre of $\Delta ^n$. Here we use the cone construction which is defined in the following way. Let $X\subseteq N^m$ be a convex  definable set and let $p\in X$. The {\it cone construction over $p$ in $X$} is a sequence of homomorphisms $z\mapsto p.z$$:\tilde{S}_*(X)\into \tilde{S}_{*+1}(X)$ defined as follows: For $n<-1$, $p.$ is defined as the trivial homomorphism and
for $n\geq -1$ and a  basis element $\sigma $, we set 
\begin{equation*}
 p.\sigma (\sum _{i=0}^{n+1}t_ie_i)=
\begin{cases} 
p \,\,\,\,\,\textmd{if} \,\,\,\,t_0=1\\
\\
\,\,\,\,\,\,\,\,\,\,\,\,\,\,\,\,\,\,\,\,\,\,\,\,\,\,\,\,\,\,\,
\\
t_0p+(1-t_0)\sigma (\sum _{i=1}^{n+1}\frac{t_i}{1-t_0}e_i)\,\,\,\, \textmd{if}\,\,\, t_0\neq 1.
\end{cases}
\end{equation*}

In the classical case we apply the Lebesgue number property to the repeated barycentric subdivision operator 
$${\rm Sd}^k=({\rm Sd}^k_n)_{n\in \ZZ}:\tilde{S}_*(X)\into \tilde{S}_*(X)$$ 
where ${\rm Sd}^k$ is the composition of ${\rm Sd}$ with itself $k$ times, to prove the following lemma. 

\begin{lem}\label{lem clas exci}
Suppose that $X$ is a topological space and let $U$ and $V$ be open subsets of $X$ such that $X=U\cup V$. If $z\in \tilde{S}_n(X)$, then there is a sufficiently large $k\in \NN$ such that ${\rm Sd}^k_n(z)\in \tilde{S}_n(U)+\tilde{S}_n(V)$.
\end{lem}

This lemma implies the excision axiom. In the o-minimal case Woerheide replaces ${\rm Sd}^k$ by the {\it subdivision operator} $${\rm Sd}_i^K:\tilde{S}_i(X)\into \tilde{S}_i(X)$$ where $(\Phi, K)$ is a definable triangulation of $X$. The subdivision operator is defined by
$${\rm Sd}_i^K(z)=({\rm Sd}z)_{\sharp }(\gamma _{i}^{n})_{\sharp }(\Phi ^{-1})_{\sharp
}\tau _{K}F_n\langle e_{n-i},\dots ,e_n\rangle $$
where $F_n:\tilde{C}_*(E^n)\into \tilde{C}_*(K)$ is the o-minimal simplicial chain map induced by $\Phi :E^n\into K$, $\gamma ^n_i:\Delta ^n\into \Delta ^i$ is defined by 
$$\gamma ^n_i(\sum _{j=0}^na_je_j)=\sum _{j=0}^i(a_{n-i+j}+\frac{\sum
_{k=0}^{n-i-1}a_k}{i+1})e_j$$
and $E^n$ is the standard simplicial complex such that $|E^n|=\Delta ^n$.

Woerheide proves the following lemma which, as in the classical case, implies the o-minimal excision axiom.

\begin{lem}[\cite{Wo}]\label{lem omin exci}
Suppose that $X$ is a definable set and let $U$ and $V$ be open  definable subsets of $X$ such that $X=U\cup V$. If $z\in \tilde{S}_n(X)$, then there is a definable triangulation $(\Phi ,K)$ of 
$\Delta ^n$ compatible with $E^n$ such that ${\rm Sd}^K_n(z)\in \tilde{S}_n(U)+\tilde{S}_n(V)$.
\end{lem}

Woerheide's construction easily gives, as in the classical case (\cite{d} Chapter VI, Section 7), o-minimal singular homology with 
coefficients in ${\QQ}$. Indeed, if $f:X\into Y$ is a definable 
continuous map, one defines {\it o-minimal singular homology with 
coefficients in ${\QQ}$} by
$$H_m(X;{\QQ})=H_m(S_*(X)\otimes {\QQ})$$
and $f_*:H_m(X;{\QQ})\into H_m(Y;{\QQ})$ is the homomorphism induced by 
$f_{\sharp }\otimes {\rm id}$. This o-minimal homology with coefficients 
in ${\QQ}$ satisfies the corresponding Eilenberg-Steenrod axioms. 
We often apply the Universal Coefficient 
theorem and identify $H_m(X)\otimes {\QQ}$ with $H_m(X;{\QQ})$ (and the 
corresponding $f_*$'s) as ${\QQ}$-vector spaces.

Similarly, as in the classical case (\cite{d} Chapter VI, Section 7), we have the {\it o-minimal singular cohomology with coefficients in ${\QQ}$}
$$H^m(X;{\QQ})=H_{-m}({\rm Hom}(S_*(X),{\QQ}))$$
with homomorphism  $f^*:H^m(Y;{\QQ})\into H^m(X;{\QQ})$  induced by 
${\rm Hom}(f_{\sharp }, {\QQ})$. This o-minimal cohomology with coefficients in ${\QQ}$ 
satisfies the corresponding corresponding Eilenberg-Steenrod axioms. 
We often apply the Universal Coefficients theorem and identify ${\rm Hom}(H_m(X),{\QQ})$ with
$H^m(X;{\QQ})$ (and the corresponding $f^*$'s) as ${\QQ}$-vector spaces.


By construction of $(H_*,d_*)$ and $(H^*,d^*)$ one can also develop the theory of products 
for the o-minimal singular homology and cohomology in the same purely algebraic way as in the classical case  (\cite{d} Chapter VI and VII). For completeness we  recall this  in the Appendix (Section \ref{section appendix}  below) since it will be used in the proof of our main result. 

For further details on o-minimal singular homology the reader should see the paper \cite{ew} by the authors.

\end{section}

\begin{section}{O-minimal relative Poincar\'e duality}\label{section omin rel pd}

Before we prove the o-minimal relative Poincar\'e duality we introduce some notation and recall orientation theory for definable manifolds.

In this paper, by a definable manifold we always mean an affine Hausdorff definable manifold, i.e., a definable subset $X$ of $N^k$ with a cover  by relatively open definable subsets $U_1,\dots ,U_l$ such that, for each $i=1,\dots ,l$, there is a definable homeomorphism $\phi _i:U_i\into V_i$ where $V_i$ is an open definable subset of $N^n$ and, for all $j=1,\dots ,l$, the map $\phi _i\circ \phi _j^{-1}:\phi _j(U_i\cap U_j)\into \phi _i(U_i\cap U_j)$ is a definable homeomorphism. A definable manifold (or a definable set) $X$ is definably compact  if it is a closed and bounded subset of $N^k$ (see \cite{ps}) and $X$ is definably connected if and only if it is not the union of two disjoint clopen definable subsets. 

Let $X$ be a definable manifold of dimension $n$. We call a finite collection $(W_l,h_l)_{l\in L}$ of open definable subsets $W_l$ of $X$ together with the definable homeomorphisms $h_l:W_l\into B_n(0,\epsilon _l)\subseteq R^n$  {\it definable charts of $X$ by open balls}. In this context it is natural to call each $W_l$ a {\it definable sub-ball of $X$} and a definable subset $U$ of $X$ of the form $h_l^{-1}(B_n(0,\delta ))$ with $0<\delta <\epsilon _l$ a {\it definable proper sub-ball of $W_l$ (or of X))} since we will have  a definable homeomorphism from  the closure $\bar{U}$ of $U$ in $X$ into the closed unit ball in $R^n$ sending $\bar{U}-U$ into the unit $(n-1)$-sphere.
 
 In this context we have the following fundamental result:
 
 \begin{thm}[\cite{e3}, \cite{bo}, \cite{w}]\label{thm cov}
If $X$ is a definable manifold of dimension $n$, then $X$ can be covered by finitely many definable sub-balls of  $X$. In particular, if $A\subseteq X$ is a definably compact definable subset of $X$, then $A$ can be covered by finitely many definable proper sub-balls of $X$.
\end{thm}

This theorem is used to define orientation theory for definable manifolds:

\begin{defn}[\cite{bo}, \cite{beo}]\label{defn orientation new}
{\em
An {\it orientation} on a definable manifold $X$ of dimension $n$  is a map 
$$s:X\into \sqcup _{x\in X}H_n(X,X-x;{\ZZ})$$ 
which assigns to each $x\in X$ a generator $s(x)$ of $H_n(X,X-x;{\ZZ})\simeq {\ZZ}$ and is such that for every definable proper sub-ball $U$ of $X$ there is a   class $\alpha _U\in H_n(X,X-U;{\ZZ})$ such that $s(u)=j^{U}_u(\alpha _U)$ for each $u\in U$, where $j^{U}_u:H_n(X,X-U;{\ZZ})\into H_n(X,X-u;{\ZZ})$ is the homomorphism induced by the inclusion.  
}
\end{defn}

Theorem \ref{thm cov} is used together with classical arguments to prove:

\begin{thm}[\cite{bo}]\label{orient fund class}
Suppose that $X$ is a definable manifold of dimension $n$ with an orientation $s$. If $A$ is a definably compact definable subset of $X$, then $s_{|A}$ is uniquely determined by a fundamental class $\zeta _A$ in $H_n(X,X-A;{\ZZ})$ as $s_{|A}(x)=j^A_x(\zeta _A )$ where $j^A_x:H_n(X,X-A;{\ZZ})\into H_n(X,X-x;{\ZZ})$ is the homomorphism induced by the inclusion. 
\end{thm}


We now proceed towards the proof of the o-minimal relative Poincar\'e duality. But we will require the following lemma.

\begin{lem}\label{lem limit of cohomo}
Let $X$ be a definable manifold. If $L\subseteq K$ are definably compact definable subset 
of $X$, then $H^*(K,L;{\QQ})$ is isomorphic to the direct limit 
$\lim _{(U,V)\in \Omega (K,L)}H^*(U,V;{\QQ})$ where $\Omega (K,L)$ is the set of pairs 
$(U,V)$ such that $U$ (resp., $V$) is an open definable neighbourhood of $K$ (resp., $L$) 
in $X$ directed by reversed inclusion.
\end{lem}

\pf
We first show that if $K$ is a definably compact definable subset of $X$, then 
$H^*(K;{\QQ})$ is isomorphic to the direct limit 
$\lim _{(X,V)\in \Omega (X,K)}H^*(V;{\QQ})$. 

In fact, if $V$ is an open definable neighbourhood of $K$ in $X$, by \cite{vdd} Chapter 
VIII, 3.3 and 3.4, there is an open definable neigbourhood $U$ of $K$ in $X$ such that 
$U\subseteq V$ such that $K$ is a definable deformation retract of $U$. Hence, the 
inclusion $K\into U$ induces an isomorphism $H^*(U;{\QQ})\into H^*(K;{\QQ})$. 

The general case stated in the lemma follows from the special case together with the 
exactness axiom.
\qed



\begin{thm}\label{thm poincare duality}
Assume that $X$ is an orientable definable manifold of dimension $n$ 
and let 
$B\subseteq A$ be definably compact definable subsets of $X$. \label{not poi disom a} 
Then, for every $q\in \ZZ$, there is an isomorphism 
$$D_{X,A}:H^q(A,B;{\QQ})\into H_{n-q}(X-B,X-A;{\QQ})$$
which is natural with respect to inclusions of pairs of definably 
compact definable subsets of $X$. 
\end{thm}

\pf
First we observe that if $K_1,K_2\subseteq X$ are definably compact definable subsets and 
the theorem holds for $(K_1,\emptyset ), (K_2,\emptyset )$ and $(K_1\cap K_2,\emptyset )$, 
then the theorem holds for $(K_1\cup K_2,\emptyset )$. 

\medskip
For $i=1,2$, let $V_i$ be a definable open neighbourhood of $K_i$ in $X$. Consider the 
diagram
{\tiny
\[\xymatrix{
H^{q-1}(V_1;{\QQ})\oplus H^{q-1}(V_2;{\QQ})\ar[d]^{\cap \zeta _{K_1}\oplus \zeta _{K_2}}\ar[rr]&&
H^{q-1}(V;{\QQ})\ar[d]^{\cap \zeta _K}\ar[rr]^{d^*}&&
H^{q}(W;{\QQ})\ar[d]^{\cap \zeta _L}   \\
H_{n-q+1}(K'_1;{\QQ})\oplus H_{n-q+1}(K'_2;{\QQ})\ar[rr]&&
H_{n-q+1}(K';{\QQ})\ar[rr]^{d_* }&&
H_{n-q}(L';{\QQ})
.}
\]
}
\noindent
where $K=K_1\cap K_2$, $L=K_1\cup K_2$, $V=V_1\cap V_2$, $W=V_1\cup V_2$ and for 
$A\in \{K_1,K_2,K,L\}$ we use $A'$ to denote the pair $(X,X-A)$.

In this diagram the rows are from the Mayer-Vietoris sequence (\cite{Wo} or \cite{e1}) 
and therefore are exact. The first and the third squares are commutative by naturality of 
cap product (Theorem \ref{thm cap product} (1)). By excision, $H_{n-q+1}(K';{\QQ})\simeq H_{n-q+1}(V,V-K;{\QQ})$ 
and $H_{n-q}(L';{\QQ})\simeq H_{n-q}(W,W-L;{\QQ})$. Hence, by Proposition \ref{prop extra cap product}
 taking $X=W$, $X_i=V_i$, $Y_i=X-K_i$ and $\alpha =\zeta _{L}$, 
we see that the second square in this diagram is commutative.

The mapping from $\Omega (K_1,K_2)$ into $\Omega (X,K_i)$ (resp., $\Omega (X,K)$ and 
$\Omega (X,L)$) which sends $(V_1,V_2)$ to $(X,V_i)$ (resp., $(X,V_1\cap V_2)$ and 
$(X,V_1\cup V_2)$) is cofinal. If we pass to the limit, by Lemma \ref{lem limit of cohomo}, 
we get the diagram 
{\tiny
\[\xymatrix{
H^{q-1}(K_1;{\QQ})\oplus H^{q-1}(K_2;{\QQ})\ar[d]^{D_{X,K_1}\oplus D_{X,K_2}}\ar[r]&
H^{q-1}(K;{\QQ})\ar[d]^{D_{X,K}}\ar[r]^{d^*}&
H^{q}(L;{\QQ})\ar[d]^{D_{X,L}}  \\
H_{n-q+1}(K'_1;{\QQ})\oplus H_{n-q+1}(K'_2;{\QQ})\ar[r]&
H_{n-q+1}(K';{\QQ})\ar[r]^{d_* }&
H_{n-q}(L';{\QQ})
.}
\]
}

By \cite{d} Chapter VIII, 5.21 (a purely algebraic result) this diagram is still commutative with exact rows. By 
assumption and the five lemma, the arrow with $D_{X,L}$ is an isomorphism as 
required. 

\medskip
We now show that the theorem holds for pairs of the form $(K,\emptyset )$ where $K$ is a 
nonempty definably compact  definable subset of $X$. 

\medskip
Arguing as in the proof of Case 5 in the proof of \cite{bo} Theorem 5.2, we see that there is a finite family 
$\{\emptyset ,K_1,\dots ,K_l\}$ closed under intersection of definably compact definable 
subsets of $K$ such that $K=\cup \{K_i:i=1,\dots ,l\}$ and there are finitely many 
definable proper sub-balls $U_1,\dots , U_k$ in  $X$  such that for each 
$i$ there is a $j_i$ such that $K_i\subseteq U_{j_i}$. 

The theorem holds for 
$(K,\emptyset )$ by induction on $l$. The inductive step follows from what we saw at the 
beginning of the proof. So suppose that $K=K_i$ and $U=U_{j_i}$. Since we are 
interested in the limit of the homomorphisms 
$-\cap \zeta _K:H^q(V;{\QQ})\into H_{n-q}(X,X-K;{\QQ})$
with $V\in \Omega (X,K)$, and by the excision axiom this limit is the same as the limit
 $-\cap \zeta _K:H^q(V;{\QQ})\into H_{n-q}(U,U-K;{\QQ})$
with $V\in \Omega (U,K)$, by the definable triangulation theorem (\cite{vdd}), we can assume that
$U=N^n$ and $K$ is the  geometric realization of a closed simplicial complex in $N^n$. 
Furthermore,  as explained above, by
induction on the number of closed simpleces, we can assume that $K$ is 
 the geometric realization of a closed simplex in $N^n$. The argument in 
the proof of Case 1 in the proof of \cite{bo} Theorem 5.2 shows that 
$H_{n-q}(N^n,N^n-K;{\QQ})$ is 
zero except for $q=0$ in which case it is ${\QQ}$.  On the other hand, clearly 
$H^q(K;{\QQ})$ is zero except for $q=0$ in which case it is ${\QQ}$ and by definable
retration the same holds for the cohomology of elements in 
a cofinal collection ${\mathcal C}$ of open definable sets in 
$\Omega (N^n,K)$. So the homomorphisms 
$-\cap \zeta _K:H^q(V;{\QQ})\into H_{n-q}(N^n,N^n-K;{\QQ})$
are isomorphisms for all $V\in {\mathcal C}$, and hence, the limit homomorphism 
$D_{X,K}:H^q(V;{\QQ})\into H_{n-q}(X,X-K;{\QQ})$ is an isomorphism as required.

\medskip
We now prove the general case. Consider the diagram
{\tiny
\[\xymatrix{
H^{q-1}(K;{\QQ})\ar[d]^{D_{X,K}}\ar[r]&
H^{q-1}(L;{\QQ})\ar[d]^{D_{X,L}}\ar[r]^{d^*}&
H^q(K,L;{\QQ})\ar[d]^{D_{X,K}}\ar[r] &
H^q(K;{\QQ})\ar[d]^{D_{X,K}} \\
H_{n-q+1}(K';{\QQ})\ar[r]&
H_{n-q+1}(L';{\QQ})\ar[r]^{d_* }&
H_{n-q}(L'',K'';{\QQ})\ar[r]&
H_{n-q}(K';{\QQ})
.}
\]
}
\noindent
where $L'=(X,X-L)$, $K'=(X,X-K)$, $L''=X-L$ and $K''=X-K$. In this diagram, the first 
row is exact by exactness axiom, the second row is exact by Mayer-Vietoris (\cite{Wo} 
or \cite{e1}), the first and the third squares are commutative by naturality of cap 
product (Theorem \ref{thm cap product} (1)). The second square is commutative, because it is the direct limit of the corresponding squares for open definable neighbourhoods $V\subseteq V'$ of 
$L\subseteq K$, and each of these squares commutes by Corollary \ref{cor extra cap product} 
with $X=V'$, $W=X-L$ and $U=X-K$. Therefore, the 
$5$-lemma and the theorem for $(K,\emptyset )$ and $(L,\emptyset )$ implies that the 
theorem holds for $(K,L)$.
\qed

\begin{cor}\label{cor poincare duality abs}
Let $X$ be an orientable, definably compact definable manifold of dimension $n$. Then 
for all $q\in \ZZ$, the homomorphism \label{not pd abs}
$$D_X:H^q(X;{\QQ})\into H_{n-q}(X;{\QQ}),\,\,\, D_X(\sigma )= \sigma \cap \zeta _X$$ 
is an isomorphism and determines a dual pairing 
$$\langle -,-\rangle :H^q(X;{\QQ})\otimes   H^{n-q}(X;{\QQ})\into {\QQ}$$ 
given by $\langle x,y\rangle =(x\cup y,\zeta _X)$. 
\end{cor}

\pf
The fact that $D_{X,X}(\sigma )=D_X(\sigma )= \sigma \cap \zeta _X$ is an isomorphism is a 
consequence of the proof of Theorem \ref{thm poincare duality}. Since 
$(x\cup y,\zeta _X)=(-1)^{{\rm deg}x{\rm deg}y}(y,x\cap \zeta _X)$, the Kronecker 
product $(\,\,,\,\,)$ is a dual pairing and $-\cap \zeta _X$ is an isomorphism, it 
follows that $\langle -,-\rangle $ is a dual pairing.
\qed

Another consequence of Theorem \ref{thm poincare duality} is the theory of o-minimal 
homology and cohomology transfers of continuous definable maps which we now present as 
it will be requred later.

\begin{cor}\label{thm cohomo  and homo transfers}
Suppose that $f:X\into Y$ is a continuous definable map of orientable, definably compact 
definable manifolds of dimensions $n$ and $m$ respectively.  Then there is a homomorphism
$$f^{!}:H^q(X;{\QQ})\into H^{m-n+q}(Y;{\QQ}),$$
called cohomology tranfer, which is given by $D_{Y}^{-1}\circ f_*\circ D_{X}$, and the 
following hold: (1) $(g\circ f)^{!}=g^{!}\circ f^{!}$; (2) $1_{X}^{!}={\rm id}$; (3) 
$f^{!}(f^*\alpha \cup \beta )=\alpha \cup f^{!}\beta .$ 

Similarly, there is a homomorphism
$$f_{!}:H_q(Y;{\QQ})\into H_{m-n+q}(X;{\QQ}),$$
called homology tranfer, which is given by $D_{X}\circ f^*\circ D_{Y}^{-1}$, and the 
following hold: (1) $(g\circ f)_{!}=f_{!}\circ g_{!}$; (2) $1_{X\,!}={\rm id}$; (3) 
$f_{!}(\alpha \cap \beta )=f^*\alpha \cap f_{!}\beta $; (4) 
$f_*(\alpha \cap f_!\beta )=(-1)^{(m-{\rm deg}\beta )(m-n)}f^!\alpha \cap \beta $. 
\end{cor}

\pf
This follows easily from the definitions. For details compare with \cite{d} Chapter VIII, 
Exercise 10.14 (4).
\qed
 
The remarks that follow below are also easy consequences of the definitions together 
with the properties of o-minimal singular (co)homology products.

\begin{nrmk}\label{nrmk cohomo  and homo transfers multi}
{\em
Suppose that $f:X\into Y$ and $g:Z\into W$ are continuous definable maps of orientable, 
definably compact  definable manifolds of dimensions $n, m, l$ and $k$ respectively.  
Then 
$$(f\times g)^!(\alpha \times \beta )=
(-1)^{(n+m){\rm deg}\beta +m(k-l)}f^!(\alpha )\times g^!(\beta )$$
and
$$(f\times g)_!(\sigma \times \tau )=
(-1)^{(n+m)(k-{\rm deg}\tau )}f_!(\sigma  )\times g_!(\tau ).$$
}
\end{nrmk}

\begin{nrmk}\label{nrmk cohomo  and homo transfers}
{\em
Suppose that $f:X\into Y$ is a continuous definable map of orientable, definably compact 
definable manifolds of dimensions $n$. Then $f_*\circ f_!={\rm deg}f=f^!\circ f^*$, $f_!\circ f_*={\rm deg}f$ on the image of $f_!$ and $f^*\circ f^!={\rm deg}f$ on the image of $f^*$. For the notion of 
degree ${\rm deg}f$ of a continuous definable map see \cite{eo}.
}
\end{nrmk}

\end{section}

\begin{section}{Lefschetz coincidence theorem}
\label{section lefschetz fixed point theorem}

Once we develop the theory of the Thom,  Lefschetz and Euler classes  below, we introduce the Lefschetz coincidence number of continuous definable maps and prove in a rather classical and algebraic way the Lefschetz coincidence theorem.

\begin{subsection}{The Thom,  Lefschetz and Euler classes}
\label{subsection thom, lefschetz and 
 euler classes}


Let $Y$ be an  orientable, definable manifold of dimension $n+k$, $X$ an  orientable, definably compact  definable manifold of dimension $n$ and $z:X\into Y$ a closed definable embedding. We 
assume that $z(X)$ is  orientable with the induced  orientation. 

Let $A$ be a definably compact definable subset of $X$. By Theorem 
\ref{thm poincare duality}, for all $q\in {\ZZ}$, we have an isomorphism
$$D_{X,A}\circ z^*\circ D^{-1}_{Y,z(A)}:H_q(Y,Y-z(A);{\QQ})\into H_{q-k}(X,X-A;{\QQ}).$$
In particular, we have that 
$$H_q(Y,Y-z(A);{\QQ})=0\,\,\,\,\,\,\, \textmd{for} \,\,\,\,q<k$$ 
and 
$$H_k(Y,Y-z(A);{\QQ})\simeq H_0(X,X-A;{\QQ})\simeq {\QQ}^l$$ 
where $l$ is the number 
of definably connected components of $X$ which lie in $A$. 

\begin{defn}\label{defn thom and transversal1}
{\em 
The generators $\nu _1,\dots ,\nu _l$ of $H_k(Y,Y-z(X);{\QQ})$ are called the 
{\it transverse classes}.  If $X$ is definably connected, we denote the unique 
tranverse class by  $\nu _{Y,X}$. 

The unique class $\tau _{Y,X}\in H^k(Y,Y-z(X);{\QQ})$ such $(\tau _{Y,X},\nu _i)=1$ 
for all $i=1,\dots ,l$, is called the {\it Thom class} and its image 
$\Lambda _{Y,X}=j^*(\tau _{Y,X})\in H^k(Y;{\QQ})$, where $j:Y\into (Y,Y-z(X))$ is the 
inclusion, is called the {\it Lefschtez class}. The class 
$\chi _{Y,X}=z^*(\Lambda _{Y,X})\in H^n(X;{\QQ})$ is called the {\it Euler class}.
}
\end{defn}

\begin{expl}\label{expl thom1}
{\em
Let $X$ be an  orientable, definably compact definable manifold of dimension $n$ and 
$\Delta _X:X\into X\times X$ the diagonal map and $\Delta _X\subseteq X\times X$ the 
diagonal. The Thom class $\tau _{X\times X,X}$ is denoted by 
$\tau _{X}\in H^n(X\times X,X\times X-\Delta _X;{\QQ})$. \label{not thom cl of x1} The 
Lefschtez class $\Lambda _{X\times X,X}$ is denoted by 
$\Lambda _X\in H^n(X\times X;{\QQ})$. \label{not lefs cl of x1} The Euler class 
$\chi _{X\times X,X}$ is denoted by $\chi _X\in H^n(X;{\QQ})$. \label{not euler cl of x1}  
}
\end{expl}

As in the classical case, below we set
$$z_!:=(-1)^{k(n+k-q)}D_{X,A}\circ z^*\circ D^{-1}_{Y,z(A)}.$$ 

\begin{prop}\label{prop thom class1}
Let $Y$ be an  orientable,  definable manifold of dimension $n+k$. Suppose $X$ is an  
orientable, definably compact definable manifold of dimension $n$, $z:X\into Y$ a 
closed definable embedding and $z(X)$ is  orientable with the induced  orientation. Let $A$ be a definably compact definable subset of $X$ and $W$ an open definable subset of $Y$ such that 
$z(X)-z(A)\subseteq W\subseteq Y-z(A)$. Then 
$z_*(\zeta _{X,A})=\tau _{Y,X}\cap \zeta _{Y,z(A)}$ where $z:(X,X-A)\into (Y,W)$ is the 
inclusion.
\end{prop}

\pf
First observe that by Theorem \ref{thm cohomo and homo transfers}, we have: 

\begin{clm}\label{clm thom class1}
$z_{!}(\alpha \cap \zeta _{Y,z(A)})=(-1)^{k{\rm deg}\alpha }(z^*\alpha )\cap \zeta _{X,A}$ 
for all $\alpha \in H^*(Y,W;{\QQ})$.
\end{clm}

\pf
If $\lambda $ is the image of $\alpha $ under the composition 
$$H^*(Y,W;{\QQ})\,{\rightarrow }\,H^*(z(X),z(X)-z(A);{\QQ})\,{\rightarrow }\,
H^*(z(X);{\QQ})$$
where the last arrow is induced by the isomorphism of Lemma \ref{lem limit of cohomo}, then
$\alpha \cap \zeta _{Y,z(A)}=D_{Y,z(X)}(\lambda )$ and 
$(z^*\alpha )\cap \zeta _{X,A}=D_{X}(z^*\lambda )$ by the definition of the right 
hand sides. Since $z_{!}(D_{Y,z(X)}(\lambda ))
=(-1)^{k{\rm deg}\alpha }D_{X}(z^*\lambda )$,
the claim holds.
\qed

For $U$ an open definable subset of $Y$ such that $z(X)\subseteq W\subseteq U$ and $z(X)$ 
is closed in $U$, let $\tau _{U,X}$ and $\zeta _{U,z(A)}$ be the classes obtained from 
$\tau _{Y,X}$ and $\zeta  _{Y,z(A)}$ by excision isomorphisms. 

\begin{clm}\label{clm thom class2}
If $r:(U,V)\into (X,X-A)$ is a definable retraction, i.e., $r\circ z=1_X$, where $V$ is 
an open definable subset of $W$ such that $V\subseteq U-z(A)$, then 
$r_*(\tau _{U,X}\cap \zeta _{U,z(A)})=\zeta _{X,A}$.
\end{clm}

\pf
We start by proving the claim for $A$ a point $x$. Let $\mu \in H^n(X,X-x;{\QQ})$ be 
such that $(\mu ,\zeta _{X,x})=1$. Then $(r^*\mu )\cap \zeta _{U,z(x)}=(-1)^{kn}\nu _{U,x}$ 
where $\nu _{U,x}$ is the transverse class of $x$. Indeed, by Claim \ref{clm thom class1}, 
\begin{eqnarray*}
(-1)^{kn}z_!((r^*\mu )\cap \zeta _{U,z(x)})&=&((z^*\circ r^*)\mu )\cap \zeta _{X,x}\\
&=&\mu \cap \zeta _{X,x}.
\end{eqnarray*}
 But $\mu \cap \zeta _{X,x}$ equals $(\mu , \zeta _{X,x})=1$ 
times the homology class of $x$. Hence, by definition of transverse classes, 
$(r^*\mu )\cap \zeta _{U,z(x)}=(-1)^{kn}\nu _{U,x}$ as required. 

We have 
\begin{eqnarray*}
(\mu ,r_*(\tau _{U,X}\cap \zeta _{U,z(x)}))&=&(r^*\mu ,\tau _{U,X}\cap \zeta _{U,z(x)})\\
&=&(r^*\mu \cup \tau _{U,X}, \zeta _{U,z(x)})\\ 
&=&(-1)^{kn}(\tau _{U,X}, (r^*\mu )\cap \zeta _{U,z(x)})\\
&=&(\tau _{U,X}, \nu _{U,x})\\
&=&1. 
\end{eqnarray*}

Thus, since $r_*(\tau _{U,X}\cap \zeta _{U,z(x)})$ is 
a multiple of $\zeta _{X,x}$, this proves that $r_*(\tau _{U,X}\cap \zeta _{U,x})=
\zeta _{X,x}$.

For the general case, let $x\in A$ and consider the inclusions 
$l:(U,V)\into (U,r^{-1}(z(X)-z(x)))$ and $i:(X,X-A)\into (X,X-x)$. Then we have a 
commutative diagram 
\[\xymatrix{
H_*(U,V;{\QQ})\ar[d]^{r_*}\ar[rr]^{l_*}&&
H_*(U,r^{-1}(X-x);{\QQ})\ar[d]^{r_*}  \\
H_*(X,X-A;{\QQ})\ar[rr]^{i_*}&&
H_*(X,X-x;{\QQ})
.}
\] 

Then 
\begin{eqnarray*}
i_*(r_*(\tau _{U,X}\cap \zeta _{U,z(A)}))&=&r_*\circ l_*(\tau _{U,X}\cap \zeta _{U,z(A)})\\
&=&r_*(\tau _{U,X}\cap \zeta _{U,z(x)})\\
&=&\zeta _{X,x} 
\end{eqnarray*}
using naturality of cap products and the first part of the proof. But 
then, by definition of $\zeta _{X,A}$, we have $r_*(\tau _{U,X}\cap \zeta _{U,z(A)})=
\zeta _{X,A}$.
\qed

Since $z(X)$ is closed in $Y$, by \cite{vdd} Chapter VIII, 3.3, there is a definable 
retraction $r:U\into X$ where $U$ is an  open definable subset of $Y$ such that 
$z(X)\subseteq U$ and $z(X)$ is closed in $U$. Let $V$ be a definable neighbourhood of 
$z(X)-z(A)$ in $r^{-1}(X-A)\cap W$. Then by \cite{vdd} Chapter VIII, 3.4, the composition 
$(U,V)\stackrel{r}{\rightarrow }(X,X-A)\stackrel{z}{\rightarrow }(Y,W)$ is definably 
homotopic to the inclusion $i:(U,V)\into (Y,W)$. Hence, $i_*=z_*\circ r_*$ and, by 
Claim \ref{clm thom class2} and naturality of cap products, 
\begin{eqnarray*}
z_*(\zeta _{X,A})&=&z_*\circ r_*(\tau _{U,X}\cap \zeta _{U,z(A)})\\
&=&i_*(i^*\tau _{Y,X}\cap \zeta _{U,z(A)})\\
&=&\tau _{Y,X}\cap i_*\zeta _{U,z(A)}\\
&=&\tau _{Y,X}\cap \zeta _{Y,z(A)}
\end{eqnarray*} 
as required.
\qed

The proof of our next result is purely algebraic using Proposition \ref{prop thom class1}.

\begin{prop}\label{prop thom euler class}
Suppose that $X$ is an  orientable, definably compact definable manifold of dimension 
$n$. Let $\{b_i:i\in I\}$ be a basis of $H^*(X;{\QQ})$ and $\{\hat{b}_i:i\in I\}$ the 
dual basis of $H^*(X;{\QQ})$, i.e., $\langle \hat{b}_i,b_j \rangle  =\delta _{ij}$ for 
all $i,j\in I$. Then
$$\Lambda _X=\sum_{i\in I}(-1)^{{\rm deg}b_i}\hat{b}_i\times b_i\,\,{\rm and}\,\,\chi _X=
\sum_{i\in I}(-1)^{{\rm deg}b_i}\hat{b}_i\cup b_i.$$

Furthermore, $(\chi _X,\zeta _X)=\chi (X)$, the o-minimal Euler-Poincar\'e characteristic 
of $X$, and, $\Delta _{X*}(\zeta  _X)=\Lambda _X\cap \zeta _{X\times X}$.
\end{prop}

\pf
By Proposition \ref{prop thom class1} and naturality of cap products, we have 
\begin{eqnarray*}
\Delta _{X*}(\zeta _X)&=&\tau _X\cap \zeta _{X\times X,\Delta _X}\\
&=&\tau _X\cap j_*(\zeta _{X\times X})\\
&=&j^*(\tau _X)\cap \zeta _{X\times X}.
\end{eqnarray*} 
Thus $\Delta _{X*}(\zeta _X)=\Lambda _X\cap \zeta _{X\times X}$.

We start by proving the following claim, where we are using here the K\"unneth formula 
for o-minimal singular cohomology to express elements of $H^*(X\times X;{\QQ})$.

\begin{clm}\label{clm thom euler class}
Suppose that $\sigma =\sum _{l,k\in I}A_{l,k}\hat{b}_l\times b_k$ is an element of 
$H^*(X\times X;{\QQ})$. Then  $\langle b_i\times \hat{b}_j, \sigma  \rangle  $ 
equals $(-1)^{-n{\rm deg}b_i}A_{i,j}$.
\end{clm}

\pf
We have that $\langle b_i\times \hat{b}_j, \sigma  \rangle  $ equals 
$\sum A_{l,k}\langle b_i\times \hat{b}_j, \hat{b}_l\times b_k \rangle  $. But 
\begin{eqnarray*} 
\langle b_i\times \hat{b}_j, \hat{b}_l\times b_k \rangle &=&((b_i\times \hat{b}_j)\cup (\hat{b}_l\times b_k), \zeta _{X\times X}) \\
&=&(-1)^{(n-{\rm deg}b_j)(n-{\rm deg}b_l)}((b_i\cup \hat{b}_l)\times 
(\hat{b}_j\cup b_k), \zeta _{X\times X})\\
&=& ((b_i\cup \hat{b}_l)\times (\hat{b}_j\cup b_k), \zeta _{X\times X})\\
&=&(-1)^{n((n-{\rm deg}b_j)+{\rm deg}b_k)}(b_i\cup \hat{b}_l,\zeta _X)(\hat{b}_j\cup b_k,\zeta _X).
\end{eqnarray*}
(Where in these equalities we used: definition of duality pairing, multiplicativity of cup and cross products, $\zeta _{X\times X}=\zeta _X\times \zeta _X$ and duality of cross products respectively.)
 Finally, $(b_i\cup \hat{b}_l,\zeta _X)=
(-1)^{{\rm deg}b_i(n-{\rm deg}b_l)}(\hat{b}_l\cup b_i,\zeta _X)=
(-1)^{{\rm deg}b_i(n-{\rm deg}b_l)}\delta _{l,i}$ and $(\hat{b}_j\cup b_k,\zeta _X)=
\delta _{j,k}$. Putting all of this together and using the fact that ${\rm deg}b_l=
{\rm deg}b_k$, we see that $\langle b_i\times \hat{b}_j, \sigma  \rangle  $ equals 
$(-1)^{-n{\rm deg}b_i}A_{i,j}$.
\qed

Suppose that $\Lambda _X=\sum _{l,k\in I}A_{l,k}\hat{b}_l\times b_k$. We are going to 
compute $\langle b_i\times \hat{b}_j, \Lambda _X \rangle  $ in two ways. By Claim 
\ref{clm thom euler class}, $\langle b_i\times \hat{b}_j, \Lambda _X \rangle  =
(-1)^{-n{\rm deg}b_i}A_{i,j}$. 

On the other hand, by definition, $\langle b_i\times \hat{b}_j, \Lambda _X \rangle  $ 
is $((b_i\times \hat{b}_j)\cup \Lambda _X,\zeta _{X\times X})$ which equals 
$(b_i\times \hat{b}_j,\Lambda _X\cap \zeta _{X\times X})$. Using the definition of 
$\Lambda _X$, the last expression is equal to 
$(b_i\times \hat{b}_j,\Delta _{X*}(\zeta _{X}))$. By naturality of Kronecker product, 
this is $(\Delta _X^*(b_i\times \hat{b}_j),\zeta _{X})=(b_i\cup \hat{b}_j,\zeta _{X})=
(-1)^{{\rm deg}b_i(n-{\rm deg}b_j)}(\hat{b}_j\cup b_i,\zeta _{X})$ by the relation 
between cup and cross product and the skew commutativity of cup products. But by 
definition the last expression is 
$(-1)^{{\rm deg}b_i(n-{\rm deg}b_j)}\langle \hat{b}_j,b_i \rangle  =
(-1)^{{\rm deg}b_i(n-{\rm deg}b_j)}\delta _{ij}$. Thus, 
$A_{i,j}=(-1)^{{\rm deg}b_i}\delta _{i,j}$ as required.
 
Since $\chi _X=\Delta ^*_X(\Lambda _X)$, the description of $\chi _X$ follows from the 
relation between cup and cross product. Also, $(\chi _X,\zeta _X)=
\sum (-1)^{{\rm deg}b_i}(\hat{b}_i\cup b_i,\zeta _X)=
\sum (-1)^{{\rm deg}b_i}\delta _{i,j}=\chi (X)$. 
\qed

\end{subsection}

\begin{subsection}{Lefschetz coincidence theorem}
\label{subsection lefschetz fixed point theorem}


In this subsection, $X,Y$ and $Z$ will be definably connected, definably compact,  
orientable definable manifolds of dimension $n,m$ and $k$ respectively. Also, 
$\Delta _X:X\into X\times X$ will denote the natural inclusion of $X$ into its 
diagonal $\Delta _X$.

\begin{defn}\label{defn lefschetz isomorphism}
{\em
Let $L^p(X;{\QQ})={\rm Hom}  (H^{n-p}(X;{\QQ}),H^{n-p}(X;{\QQ}))$ and let 
\label{not hom rg of cohomo}
$$L^*(X;{\QQ})=\sum _{p=0}^nL^p(X;{\QQ}).$$ 
For each $p$, let $\{b_i:i\in I_p\}$ be a basis of $H^{n-p}(X;{\QQ})$ and let 
$\{\hat{b}_i:i\in I_p\}$ be the dual basis on $H^p(X;{\QQ})$. Then we have a canonical 
isomorphism 
$$k^p:H^p(X;{\QQ})\otimes   H_{p}(X;{\QQ})\into L^p(X;{\QQ})$$ which sends 
$\sum _{i,j\in I_p}A_{i,j}\hat{b}_i\otimes   D_X(b_j)$ into the element of 
$L^p(X;{\QQ})$ whose matrix relative to the fixed basis is $(A_{i,j})_{i,j\in I_p}$. 
The isomorphisms $k^p$ induce a canonical isomorphism
$$k:\sum _{p=0}^nH^p(X;{\QQ})\otimes   H_{p}(X;{\QQ})\into L^*(X;{\QQ})$$ given by 
$k=\sum _{p=0}^n(-1)^pk^p$.
The {\it Lefschetz isomorphism for $X$} is the isomorphism of  ${\QQ}$-modules 
$$ \lambda _X:L^*(X;{\QQ})\into H^n(X\times X;{\QQ})$$ \label{not lefs isom} given by 
$\lambda _X=
\alpha '\circ (1^*_X\otimes   D^{-1}_X)\circ
k^{-1}$ where $\alpha '$ is the K\"unneth isomorphism for o-minimal singular cohomology 
and $D^{-1}_X$ is the inverse of the Poincar\'e duality isomorphism (see Theorem 
\ref{thm poincare duality}). 
}
\end{defn}

\begin{nrmk}\label{nrmk lefschtez isom}
{\em
Note that by Proposition \ref{prop thom euler class} and the definition of $\lambda _X$, 
we have $\Lambda _X=\lambda _X(1^*_X)$.
}
\end{nrmk}

\begin{lem}\label{lem trace and lefschetz isomorphism}
Let ${\rm Tr}:L^*(X;{\QQ})\into {\QQ}$ \label{not graded trace map} be the linear map 
given by ${\rm Tr} \sigma =\sum
_{p=0}^n(-1)^p{\rm tr}_p \sigma ^p$ where $\sigma =\sum _{p=0}^n \sigma ^p$, 
$\sigma ^p\in L^p(X;{\QQ})$. Then 
$${\rm Tr}\sigma =( \Delta _X^*\lambda _X(\sigma ),\zeta _X).$$
\end{lem}

\pf
It is enough to consider $\sigma =k(\beta \otimes   D_X\gamma )$ with 
$\beta \in H^p(X;{\QQ})$ and $\gamma \in H^{n-p}(X;{\QQ})$. Then, by ordinary linear 
algebra
\begin{eqnarray*}
{\rm Tr}\sigma &=&(-1)^{p+p}(\beta, D_X\gamma )\\
&=&(\beta ,D_X\gamma )\\
&=&(\beta ,\gamma \cap \zeta _X)\\
&=& (\beta \cup \gamma ,\zeta _X)\\
&=&(\Delta _X^*\alpha '(\beta \otimes   \gamma ),\zeta _X)\\
&=&(\Delta _X^*\lambda _X(\sigma ),\zeta _X).
\end{eqnarray*}
\qed

\begin{lem}\label{lem lambda and functions}
Let $\sigma \in L^*(Y;{\QQ})$ and let $f,g:X\into Y$ be continuous definable maps where
${\dim }X={\dim }Y$. Then
$$(f\times g)^*(\lambda _Y(\sigma ))=\lambda _X(f^*\circ \sigma \circ g^!).$$
\end{lem}

\pf
It is enough to take $\sigma =k(\alpha \otimes   D_Y\beta )$ with 
$\alpha \in H^p(Y;{\QQ})$ and $\beta \in H^{n-p}(Y;{\QQ})$. We have 
\begin{eqnarray*} 
(f^*\circ \sigma \circ g^!)(\gamma )&=&(-1)^{p}(g^!(\gamma ), D_Y\beta )f^*(\alpha )\\
&=&(-1)^{p}(\gamma ,D_Xg^*(\beta))f^*(\alpha)\\
&=&[k(f^*\alpha \otimes   D_Xg^*\beta)](\gamma ).
\end{eqnarray*}
for all 
$\gamma \in H^p(X;{\QQ})$ . Therefore, $\lambda _X (f^*\circ \sigma \circ g^!)=
(f\times g)^*\circ \alpha '(\alpha \otimes   \beta)=(f\times g)^*(\lambda _Y
(\sigma )).$
\qed

\begin{defn}\label{defn coincidence number}
{\em
Let $f,g:X\into Y$ be continuous definable maps and suppose that ${\dim }X={\dim }Y$. 
The {\it Lefschetz coincidence number of $f$ and $g$} is defined by \label{not lefs coin}
$$\lambda (f,g;{\QQ})=\sum _{p=0}^n(-1)^p{\rm tr}_p(f^*\circ g^!).$$ 
Note that if $X=Y$, then $\lambda (f,1_X;{\QQ})$ is denoted by $\lambda (f;{\QQ})$.
}
\end{defn}

\begin{nrmk}\label{nrmk coincidence number}
{\em
Remark \ref{nrmk lefschtez isom} and Lemmas \ref{lem trace and lefschetz isomorphism}, 
\ref{lem lambda and functions} imply that 
$$\lambda (f,g;{\QQ})=(\Delta _X^*\circ (f\times g)^*(\Lambda _Y),\zeta _X).$$ 
Thus  $\lambda (f,g;{\QQ})=(-1)^n\lambda (g,f;{\QQ})$.  Since 
${\rm tr}_p(AB)={\rm tr}_p(BA)$, we also have 
$\lambda (f,g;{\QQ})=\sum _{p=0}^n(-1)^p{\rm tr}_p(f^!\circ g^*)$. Clearly, 
$f^*\circ g^!=D_X^{-1}\circ f_!\circ g_*\circ D_X$. So, 
$\lambda (f,g;{\QQ})=\sum _{p=0}^n(-1)^p{\rm tr}_p(f_!\circ g_*)$. Similarly, 
$g^!\circ f^*=D_Y^{-1}\circ f_!\circ g_*\circ D_Y$ and so 
$\lambda (f,g;{\QQ})=\sum _{p=0}^n(-1)^p{\rm tr}_p(f_*\circ g_!)$.

If $h:Z\into X$ is a third continuous definable map and ${\rm dim}Z={\rm dim}X$, 
then by Remark \ref{nrmk cohomo  and homo transfers}, 
$\lambda (f\circ h, g\circ h;{\QQ})=({\rm deg}h)\lambda (f,g;{\QQ})$. In particular, 
$\lambda (f,f;{\QQ})=({\rm deg}f)E(X)$ (where $E(X)$ is the o-minimal Euler characteristic of $X$, see \cite{vdd} and \cite{bo}). 
}
\end{nrmk}

We are now ready to prove the main theorem of the paper.

\vspace{.15in}
\noindent {\bf Proof of Theorem \ref{thm lct}}: We have 
$\lambda (f,g;{\QQ})=(\Delta _X^*\circ (f\times g)^*(\Lambda _Y),\zeta _X).$ If 
there is no $x\in X$ such that $f(x)=g(x)$, then we have a factorisation
\[\xymatrix{
X\ar[d]^{\Delta _{X}}\ar[rr]^{f\times g}&&
Y\times Y\\
X\times X\ar[rr]^{f\times g}&&
Y\times Y-\Delta _Y\ar[u]^{i}
.}
\]
where $i$ is the inclusion. 
Since $H^n(i)H^n(j)=0$ and $\Lambda _Y=H^n(j)(\tau_Y)$, we have 
$0=
\Delta _X^*\circ (f\times g)^* \circ i^*(\Lambda _Y)$
$=(f\times g)^*(\Lambda _Y)$ and therefore $\lambda (f,g;{\QQ})=0$.
\qed


\vspace{.15in}
We end the subsection with another characterization of the Lefschetz coincidence number and 
yet another classical proof of Theorem \ref{thm lct}.

\begin{prop}\label{prop omin inter}
Let $X$ be an  orientable, definably compact definable manifold of dimension $n$. Then 
there is a graded bilinear map (called the intersection product)
$$\cdot :H_p(X;{\QQ})\otimes   H_q(X;{\QQ})\into H_{p+q-n}(X;{\QQ})$$ defined by
$$\alpha \cdot \beta =D_X(D_X^{-1}(\alpha )\cup D_X^{-1}(\beta ))$$
and such that the following hold: 
\begin{enumerate}
\item {\bf Naturality.} $f_!(\alpha \cdot \beta )=
f_!(\alpha )\cdot f_!(\beta ).$  
\item {\bf Skew commutativity.} $\alpha \cdot \beta =(-1)^{(n-{\rm deg}\alpha )(n-{\rm deg}\beta )}\beta \cdot \alpha .$
\item {\bf Associativity.} $\alpha \cdot (\beta \cdot \gamma )=
(\alpha \cdot \beta )\cdot \gamma .$  
\item {\bf Units.} $\zeta _X\cdot \beta =\beta =\beta \cdot \zeta _X.$ 
\item {\bf Multiplicativity.} $(\alpha \cdot \beta)\times (\sigma \cdot \gamma )=
(-1)^{(n-{\rm deg}\alpha )(m-{\rm deg}\gamma )}(\alpha \times \sigma)\cdot 
(\beta \times \gamma ).$
\end{enumerate}
\end{prop}

\pf
The properties of the intersection product follow easily from the definition and the 
properties of cap and cup products.
\qed

Using the relationship between cup and cross product, it is easy to prove the following 
remark.

\begin{nrmk}\label{nrmk omin inter}
{\em
If $X$ and $Y$ are  orientable, definably compact definable manifolds of dimension $n$ 
and $m$ respectively, then the intersection product satisfies:
$$\alpha \cdot \beta =(-1)^{n(n-{\rm deg}\beta )}\Delta _{X!}(\alpha \times \beta)$$
and
$$\alpha \times \sigma =(-1)^{n(m-{\rm deg}\sigma )}(p_{X!}\alpha )\cdot (p_{Y!}\sigma )$$
where $p_X:X\times Y\into X$ and $p_Y:X\times Y\into Y$ are the projections.
}
\end{nrmk}


\begin{thm}\label{thm lefschetz omin inter}
Let $f,g:X\into Y$ be continuous definable maps between  orientable, definably compact 
definable manifolds of dimension $n$. If 
$\zeta _{f}=(1_X\times f)_*\circ \Delta _{X*}(\zeta _X)$ and 
$\zeta _{g}=(1_X\times g)_*\circ \Delta _{X*}(\zeta _X)$, then
$$\lambda (f,g;{\QQ})=\epsilon _*(\zeta _f\cdot \zeta _g)$$
where $\epsilon _*:H_0(X\times Y;{\QQ})\into {\QQ}$ is the augmentation. In particular, 
if $\lambda (f,g;{\QQ})\neq 0$, then there is $x\in X$ such that $f(x)=g(x)$.
\end{thm}

\pf
Let $\gamma _f$ and $\gamma _g$ be elements in $H^n(X\times Y;{\QQ})$ such that 
$\gamma _f\cap \zeta_{X\times Y}=\zeta _f$ and $\gamma _g\cap \zeta_{X\times Y}=\zeta _g$. 
Then 
\begin{eqnarray*} 
\epsilon _*(\zeta _f\cdot \zeta _g)&=&\epsilon _*((\gamma _g\cup \gamma _f)\cap \zeta _{X\times Y})\\
&=&\epsilon _*(\gamma _g\cap (\gamma _f\cap \zeta _{X\times Y}))\\
&=&\epsilon _*(\gamma _g\cap \zeta _f)\\
&=&(\gamma _g,\zeta _f) \,\,\,(\textmd{by definition of augmentation)}\\
&=&(D_{X\times Y}^{-1}\circ (1_X\times g)_*(\Delta _{X*}(\zeta _X)), (1_X\times f)_*(\Delta _{X*}(\zeta _X)))\\
&=&(D_{X\times Y}^{-1}\circ (1_X\times g)_*\circ D_{X\times X}(\Lambda _X), (1_X\times f)_*(\Delta _{X*}(\zeta _X)))\\
&=& ((1_X^!\times g^!)(\Lambda _X), (1_X\times f)_*(\Delta _{X*}(\zeta _X)))\,\,\,\,\,\,\,\,\,\,\textmd{(i)}\\
&=&((1_X\times f)^*\circ (1_X^!\times g^!)(\Lambda _{X}),\Delta _{X*}(\zeta _X))\,\,\,\,\,\,\,\,\,\textmd{(ii)}\\
&=&((1_X^*\circ 1_X^!\times f^*\circ g^!)(\Lambda _{X}),\Delta _{X*}(\zeta _X))\,\,\,\,\,\,\,\,\,\,\,\,\,\,\,\textmd{(iii)}\\
&=&\sum (-1)^{{\rm deg}b_i}(\hat{b}_i\times f^*\circ g^!(b_i),\Delta _{X*}(\zeta _X))\,\,\,\,\,\,\,\,\textmd{(iv)}\\
&=&\sum (-1)^{{\rm deg}b_i}(\hat{b}_i\cup f^*\circ g^!(b_i),\zeta _X)\,\,\,\,\,\,\,\,\,\,\,\,\,\,\,\,\,\,\,\,\,\,\,\,\,\textmd{(v)}\\
&=&\sum (-1)^{{\rm deg}b_i}\langle \hat{b}_i,f^*\circ g^!(b_i) \rangle \\ 
&=&\sum _{p=0}^n(-1)^p{\rm tr}_p(f^*\circ g^!)
\end{eqnarray*}
where: (i)  $D_{X\times Y}=D_X\times D_Y$ and $D_{X\times X}=D_X\times D_X$; (ii) the
naturality of the Kronecker product and Proposition \ref{prop thom euler class}; (iii) 
the  naturality of cross product and Remark \ref{nrmk cohomo  and homo transfers multi}; (iv) Proposition \ref{prop thom euler class}; (v) the naturality of the Kronecker product. 

Let $\Gamma _f$ (resp., $\Gamma _g$) be the graph of $f$ (resp., $g$). Then there are 
$\sigma _f\in H^n(X\times Y,X\times Y-\Gamma _f;{\QQ})$  and 
$\sigma _g\in H^n(X\times Y,X\times Y-\Gamma _g;{\QQ})$ such that $\gamma _f$ (resp., 
$\gamma _g$) is the image of $\sigma _f$ (resp., $\sigma _g$) by the homomorphism 
induced by inclusion. If $f$ and $g$ have no coincidence, then 
$\Gamma _f\cap \Gamma _g=\emptyset $, and so, 
$\sigma _f\cup \sigma _g\in H^{2n}(X\times Y,X\times Y;{\QQ})=0$. Therefore, by 
naturality of cup products, $\gamma _f\cup \gamma _g=0$ and 
$\lambda (f,g;{\QQ})=\epsilon _*(\zeta _f\cdot \zeta _g)=0$. 
\qed

\end{subsection}

\end{section}

\begin{section}{Appendix: O-minimal  ring (co)homology theory}\label{section appendix}

By construction of $(H_*,d_*)$ and $(H^*,d^*)$ one can also develop the theory of products 
for the o-minimal singular homology and cohomology in the same purely algebraic way as in the classical case (\cite{d} Chapter VI and VII). For completeness we include here this theory. 

First recall that if $(X,A),(X,B)$ are pairs of definable sets with $A\subseteq X$ and 
$B\subseteq X$, then we call $(X;A,B)$ a {\it definable triad}. We say 
that a definable triad $(X;A,B)$  is an {\it excisive triad (with 
respect to $(H_*,d_*)$)} if the inclusion $(A,A\cap B)\into (A\cup B,B)$ 
induces isomorphisms $H_*(A,A\cap B)\simeq H_*(A\cup B,B)$. 

Let $(X,A),(Y,B)$ be pairs of definable sets with $A\subseteq X$ and 
$B\subseteq Y$. Then we will write $(X,A)\times (Y,B)$ for 
$(X\times Y,A\times Y\cup X\times B)$. 

As we pointed out in \cite{e2}, the o-minimal version of the Eilenberg-Zilber theorem (Proposition 3.2 in \cite{eo}) gives, as in \cite{d} Chapter VI, Section 12 and Chapter VII, Section 2 respectively, the following two theorems:

\begin{thm}[K\"unneth Formula for Homology, \cite{e2}]\label{thm kunneth for hom} 
Let $(X,A)$ and $(Y,B)$ be pairs of definable sets with 
$A\subseteq X$ and $B\subseteq Y$ such that  
$(X\times Y;A\times Y,X\times B)$ is an 
excisive triad. Then, for all $n\in \ZZ$, there is an isomorphism  
$$\alpha '': \sum _{i+j=n} H_i(X,A;{\QQ})\otimes H_j(Y,B;{\QQ})\into 
H_n((X,A)\times (Y,B);{\QQ}).$$
\end{thm}

The homomorphism $\alpha ''$ from Theorem \ref{thm kunneth for hom} is 
called the {\it homology (external) cross product} 
\label{not homo ext cross prod} and $\alpha ''(a\otimes b)$ is denoted 
$a\times b$. 

\begin{thm}[\cite{e2}]\label{thm homology external cross product}
The homology cross product satisfies the following properties: 
\label{not homo cross prod}
\begin{enumerate}
\item
{\bf Naturality}. $(f\times g)_*(\alpha \times \beta )=(f_*\alpha )
\times (g_*\beta ).$ 

\item
{\bf Skew-commutativity}. $t_*(\alpha \times \beta )
=(-1)^{{\rm deg}\alpha {\rm deg}\beta }\beta \times \alpha $  
where $t:X\times Y\into Y\times X$ commutes factors. 

\item
{\bf  Associativity}. $(\alpha \times \beta )\times \gamma =\alpha 
\times (\beta \times \gamma ).$ 

\item
{\bf  Units}. $1\times \alpha =\alpha \times 1=\alpha $. 

\item
{\bf Stability}. $d_*(\alpha \times \beta )=
i_{1*}(d_*\alpha \times \beta )+
i_{2*}((-1)^{{\rm deg}\alpha }\alpha \times d_*\beta )$ where 
$\alpha \in H_i(X,A;{\QQ})$, $\beta \in H_j(Y,B;{\QQ})$ and 
$i_1:(A\times Y,A\times B)\into (A\times Y\cup X\times B,A\times B)$ and 
$i_2:(X\times B,A\times B)\into (A\times Y\cup X\times B,A\times B)$ are 
the inclusions.
\end{enumerate}
\end{thm}

By dualizing the  Eilenberg-Zilber maps from the o-minimal version of the Eilenberg-Zilber theorem (Proposition 3.2 in \cite{eo}) gives, as in \cite{d} Chapter VI, Section 10, the following:

\begin{thm}[K\"unneth Formula for Cohomology]\label{thm kunneth for cohom} 
For    pairs of definable sets $(X,A)$ and $(Y,B)$ with 
$A\subseteq X$ and $B\subseteq Y$ such that  
$(X\times Y;A\times Y,X\times B)$ is an 
excisive triad, we have that, for all $n\in \ZZ$, there is an isomorphism  
$$\alpha ': \sum _{i+j=n} H^i(X,A;{\QQ})\otimes H^j(Y,B;{\QQ})\into 
H^n((X,A)\times (Y,B);{\QQ}).$$
\end{thm}

Let $(X,A)$  be a pair of definable sets with $A\subseteq X$. The {\it Kronecker product} \label{not kron prod}
$$ (\,\, ,\,\,):H^*(X,A;{\QQ})\otimes H_*(X,A;{\QQ})\into {\QQ}$$ 
is the homomorphism $(e\otimes  {\rm id})_*\circ \alpha $ where 
{\tiny
$$\alpha :H^*(X,A;{\QQ})\otimes  H_*(X,A;{\QQ})\into H_*({\rm Hom} (S_*(X,A)\otimes {\QQ},{\QQ})\otimes (S_*(X,A)\otimes   {\QQ}))$$
} 

\noindent
is the K\"unneth homomorphism from homological algebra (see \cite{d} Chapter VI, Theorem 9.13) and 
$$(e\otimes  {\rm id}):{\rm Hom}(S_*(X,A)\otimes {\QQ},{\QQ})
\otimes (S_*(X,A)\otimes  {\QQ})\into {\QQ}\otimes  {\QQ}\simeq {\QQ}$$ 
is the evaluation chain map given by $(e\otimes  {\rm id})(\sigma \otimes (a\otimes  m))=\sigma (a)\otimes  m$.

By purely algebraic arguments, compare with \cite{d} Chapter VII, 1.8 and 1.12, we see that the Kronecker product is a dual pairing satisfying: 
$$(f^*\alpha ,\beta )=(\alpha ,f_*\beta ).$$

The homomorphism $\alpha '$ from Theorem \ref{thm kunneth for cohom} is called the {\it cohomology (external) cross product} \label{not cohomo ext cross prod}and $\alpha '(a\otimes  b)$ is denoted $a\times b$. 

\begin{thm}\label{thm cohomology external cross product}
The cohomology cross product satisfies the following properties: \label{not cohomo cross prod} 
\begin{enumerate}
\item
{\bf Naturality}. $(f\times g)^*(\alpha \times \beta )=(f^*\alpha )
\times (g^*\beta ).$ 

\item
{\bf Skew-commutativity}. $t^*(\alpha \times \beta )
=(-1)^{{\rm deg}\alpha {\rm deg}\beta }\beta \times \alpha $  
where $t:X\times Y\into Y\times X$ commutes factors. 

\item
{\bf  Associativity}. $(\alpha \times \beta )\times \gamma =\alpha 
\times (\beta \times \gamma )$. 

\item
{\bf  Units}. $1\times \alpha =\alpha \times 1=\alpha $. 

\item
{\bf  Duality}. $(\alpha \times \beta ,\sigma \times \tau )
=(-1)^{{\rm deg}\tau {\rm deg}\alpha }(\alpha ,\sigma )\otimes  (\beta ,\tau )$ where 
$\sigma \in H_*(X,A;{\QQ})$, $\tau \in H_*(Y,B;{\QQ})$, $\alpha \in H^*(X,A;{\QQ})$ 
and $\beta \in H^*(Y,B;{\QQ})$.

\item
{\bf Stability}. $d^*(\alpha \times \beta )=(i_{1*}\oplus i_{2*})^{-1}((d_*\alpha \times \beta )+((-1)^{{\rm deg}\alpha }\alpha \times d_*\beta ))$ where $\alpha \in H^i(A;{\QQ})$, $\beta \in H^j(B;{\QQ})$ and $i_1:(A\times Y,A\times B)\into (A\times Y\cup X\times B,A\times B)$ and $i_2:(X\times B,A\times B)\into (A\times Y\cup X\times B,A\times B)$ are the inclusions. 
\end{enumerate}
\end{thm}

\pf
This is obtained from the proof of Theorem \ref{thm homology external cross product} by applying the 
functor ${\rm Hom}(-,{\QQ}\otimes  {\QQ})$ and using the dual of the Eilenberg-Zilber map
from the proof of Theorem \ref{thm kunneth for cohom} (compare with \cite{d} Chapter VII, Section 7 for the details).
\qed

We now introduce the cup products. Although these are equivalent to the cross products, they are usually more convenient. 

\begin{thm}\label{thm cup product}
Suppose  $(X;A_1,A_2)$ is an excisive triad of definable sets. Then we have a canonical graded bilinear map (called cup product)\label{not cup prod} 
$$\cup :H^*(X,A_1;{\QQ})\otimes H^*(X,A_2;{\QQ})
\into H^*(X,A_1\cup A_2;{\QQ})$$ 
such that:  
\begin{enumerate}
\item
{\bf Naturality}.  $f^*(\alpha \cup \beta )=(f^*\alpha )\cup (f^*\beta )$. 

\item
{\bf Skew-commutativity}.  $\alpha \cup \beta =(-1)^{{\rm deg}\alpha {\rm deg}\beta }$ $\beta \cup \alpha $.

\item
{\bf  Associativity}.  $(\alpha \cup \beta )\cup \gamma =\alpha \cup (\beta \cup \gamma )$.

\item
{\bf  Units}. $1\cup \alpha =\alpha \cup 1=\alpha $.

\item
{\bf  Multiplicativity}. $(\alpha \times \beta)\cup (\sigma \times \tau )=(-1)^{{\rm deg}\beta {\rm deg}\sigma }(\alpha \cup \sigma )\times (\beta \cup \tau )$.

\item
{\bf Stability}.  $ d^*\circ (j^*)^{-1}(\alpha \cup i*\beta )=(d^*\alpha )\cup \beta $  where $i:(A_1,A_1\cap A_2)\into (X,A_2)$ and $j:(A_1,A_1\cap A_2)\into (A_1\cup A_2,A_2)$ are the inclusions.
\end{enumerate}
\end{thm}
 
\pf
The cup product $\cup $ is the graded bilinear map given by 
$j'^*\circ D'^*\circ \gamma _*\circ \alpha $ where: $\alpha $ is 
the map from the K\"unneth formula for cochain complexes 
(see \cite{d} Chapter VI, Theorem 9.13); $\gamma $ is the chain map from 
the proof of Theorem \ref{thm kunneth for cohom} (the dual Eilenberg-Zilber map); 
$D'={\rm Hom}(D,{\QQ})$ is the cochain map from 
{\tiny 
$${\rm Hom}  (S_*(X,A_1),{\QQ})\otimes  {\rm Hom}  (S_*(X,A_2),{\QQ})\into {\rm Hom}  (\frac{S_*(X)}{S_*(A_1)+S_*(A_2)},{\QQ})$$ 
}
\noindent
where $D$ is the natural chain map (called diagonal map) given by 
{\tiny 
$$D=\zeta \circ \Delta :\frac{S_*(X)}{S_*(A_1)+S_*(A_2)}\into S_*(X,A_1)\otimes S_*(X,A_2),$$
}
\noindent 
here $\Delta: S_*(X)\into S_*(X)\into S_*(X\times X)$ is the natural chain map 
induced by the diagonal map $X\into X\times X$ and 
$\zeta :S_*(X\times X)\into S_*(X)\otimes S_*(X)$ is a 
Eilenberg-Zilber chain map; and $j'$ is the homotopy equivalence 
{\tiny
$${\rm Hom}  (\frac{S_*(X)}{S_*(A_1)+S_*(A_2)},{\QQ})\into 
{\rm Hom}  (S_*(X,A_1\cup A_2),{\QQ})$$
}
\noindent
which exists since $(X;A_1,A_2)$ is an excisive triad of definable sets.

The properties of the cup product listed above, follow from 
corresponding properties for the Eilenberg-Zilber chain equivalence. 
Since these are purely algebraic, we refer the reader to \cite{d} 
Chapter VII, Section 8 for the details.
\qed

Similarly to  the classical case (\cite{d} Chapter VII, Section 8) we also have:

\begin{nrmk}\label{nrmk relative kunneth formula}
{\em
The cohomology cross product is related to the cup product by: 
\begin{enumerate}
\item 
$\alpha \times \beta =p^*\alpha \cup q^*\beta $ where $p:(X\times Y,A\times Y)\into (X,A)$ and $q:(X\times Y,X\times B)\into (Y,B)$ are the projections (and we assume here that $(X\times Y;A\times Y,X\times B)$ is excisive)

\item  
$\alpha \cup \beta =\Delta _X^*(\alpha \times \beta )$ where $\Delta _X:(X,A_1\cup A_2)\into (X\times X,A_1\times X\cup X\times A_2)$ is the diagonal map \label{not diag map} (and we assume here that $(X\times X;A_1\times X,X\times A_2)$ is excisive).
\end{enumerate}
}
\end{nrmk}

Theorem \ref{thm cup product} implies that $H^*(X;{\QQ})$ is 
a graded ${\QQ}$-algebra under cup product, $H^*(\,\,;{\QQ})$ is a 
functor from the category definable sets into the category of graded 
skew-commutative (associative) ${\QQ}$-algebra with unit element and 
 $H^*(X,A;{\QQ})$ is a graded $H^*(X;{\QQ})$-module 
with respect to the cup product 
$\cup :H^*(X;{\QQ})\otimes H^*(X,A;{\QQ})\into H^*(X,A;{\QQ})$. 
Moreover, the cross product $\times :H^*(X;{\QQ})\otimes H^*(Y;{\QQ})
\into H^*(X\times Y;{\QQ})$ is a homomorphism of graded 
skew-commutative associative ${\QQ}$-algebras.

Another useful product is the cap product. This is in some sense dual to
the cup product. 

\begin{thm}\label{thm cap product}
Suppose that $(X;A_1,A_2)$ is an excisive triad of definable sets. Then we have a canonical graded bilinear map (called cap product) \label{not cap prod}
$$\cap :H^p(X,A_2;{\QQ})\otimes  H_{p+q}(X,A_1\cup A_2;{\QQ})
\into H_q(X,A_1;{\QQ})$$
such that:  
\begin{enumerate}
\item 
{\bf  Naturality}. $f_*((f^*\alpha )\cap \beta )=\alpha \cap (f_*\beta )$. 

\item
{\bf Associativity}. $(\alpha \cup \beta )\cap \gamma =\alpha \cap (\beta \cap \gamma )$.  

\item
{\bf  Units}. $1\cap \alpha =\alpha \cap 1=\alpha $. 

\item
{\bf  Duality}. $(\alpha \cup \beta ,\sigma )=(\alpha ,\beta \cap \sigma )$. 

\item
{\bf Multiplicativity}. $(\alpha \times \beta)\cap (\sigma \times \tau )=(-1)^{{\rm deg}\beta {\rm deg}\sigma }(\alpha \cap \sigma )\times (\beta \cap \tau )$.

\item
{\bf Stability}. $d_*(\alpha \cap \beta )=(-1)^{{\rm deg}\alpha }(i^*(\alpha )\cap (j^*)^{-1}\circ d_*(\beta )$ and $d^*(\alpha )\cap \beta +(-1)^{{\rm deg}\alpha }i_*(\alpha \cap (j^*)^{-1}\circ d_*(\beta ))=0$ where $i:(A_1,A_1\cap A_2)\into (X,A_2)$ and $j:(A_1,A_1\cap A_2)\into (X,A_2)$ are the inclusions.
\end{enumerate}
\end{thm}

\pf
The cap product is the bilinear map given by 
$E_*\circ ({\rm id}\otimes   D)_*\circ ({\rm id} \otimes   k)_*\circ \alpha $ 
where $D$ and $\alpha $ are as before, 
{\tiny
$${\rm id}:{\rm Hom}  (S_*(X,A_2),{\QQ}) \into {\rm Hom}  (S_*(X,A_2),{\QQ})$$
}
\noindent 
is the identity map, $k$ is the chain equivalence 
{\tiny
$$k:S_*(X,A_1\cup A_2)\otimes {\QQ}\into \frac{S_*(X)}{S_*(A_1)+S_*(A_2)}\otimes {\QQ}$$
}
\noindent 
and $E$ is the natural evaluation map  
{\tiny
$${\rm Hom}(S_*(X,A_2),{\QQ})\otimes   (S_*(X,A_1)\otimes S_*(X,A_2)\otimes {\QQ})\into S_*(X,A_1)\otimes {\QQ}$$
}
\noindent
given by $E(\sigma \otimes a\otimes b\otimes m)=(-1)^{{\rm deg} \sigma {\rm deg} a}a\otimes m\otimes \sigma (b)$. 
The proofs of the properties of the cap product are simple computations as before (compare with  \cite{d} Chapter VII, Section 12 for the details).
\qed

The next result is proved as in the classical case (\cite{d} Chapter VII, Section 12, Proposition 12.20) with the following change: one replaces the use of \cite{d} Chapter III, 7.3 by its o-minimal analogue given by Lemma \ref{lem omin exci}. 

Below, we denote by $Z^n(X,A;{\QQ})$ the kernel of $\partial ^n:S^n(X,A;{\QQ})\into S^{n+1}(X,A;{\QQ})$.

\begin{prop}\label{prop extra cap product}
Suppose that $X_1,X_2,Y_1,Y_2$ are open definable subsets of a definable set $X$ such that $X_1\cup Y_1=X_2\cup Y_2=X_1\cup X_2=X$. Let $X'=X_1\cap X_2$, $Y'=Y_1\cap Y_2$, $Y=Y_1\cup Y_2$ and let
$j:(X',X'\cap Y)\into (X,Y)$ denote the inclusion. For $\alpha \in H_*(X,Y';{\QQ})$, let $\alpha '\in H_*(X',X'\cap Y;{\QQ})$ denote its image  under the composition
$$H_*(X,Y';{\QQ}){\rightarrow }H_*(X,Y;{\QQ})\stackrel{(j_*)^{-1}}{\rightarrow }H_*(X',X'\cap Y;{\QQ}).$$

Then the following diagram
\[\xymatrix{ 
H^*(X';{\QQ})\ar[d]^{\delta ^*}\ar[rr]^{j_*\circ (\,\,\,\cap \alpha ')}&&
 H_*(X,Y;{\QQ}) \ar[d]^{\delta ^*} \\
H^*(X;{\QQ})\ar[rr]^{\cap \alpha }&&
H_*(X,Y';{\QQ})
.}
\]
where $\delta _*$ and $\delta ^*$ are the Mayer-Vietoris boundaries, is commutative.
\end{prop}

\pf
By the proof of the cohomology Mayer-Vietoris sequence (\cite{Wo} or \cite{e1}), we have that $\delta ^*$ is the composition
$$H^*(X';{\QQ})\stackrel{d^*}{\rightarrow }H^*(X_1,X';{\QQ})\stackrel{(l^*)^{-1}}{\rightarrow }H^*(X,X_2;{\QQ}){\rightarrow }H^*(X;{\QQ})$$ where $l:(X_1,X')\into (X,X_2)$ is the inclusion. Similarly, by the proof of the homology Mayer-Vietoris sequence (\cite{Wo} or \cite{e1}), we have that $\delta _*$ is the composition
$$H_*(X',Y;{\QQ})\stackrel{c_*\circ d_*}{\rightarrow }H_*(Y,Y_1;{\QQ})\stackrel{(k_*)^{-1}}{\rightarrow }H_*(Y_2,Y';{\QQ}){\rightarrow }H_*(X,Y';{\QQ})$$ where $a:(Y,\emptyset )\into (Y,Y_1)$ and $k:(Y_2,Y')\into (Y,Y_1)$ are the inclusions.

Let $b:(Y_2,X_2\cap Y_2)\into (X,X_2)$ be the inclusion and let $\alpha _1\in H_*(Y_2,Y'\cup (X_2\cap Y_2);{\QQ})$ denote the image of $\alpha $ under the composition
$$H_*(X,Y';{\QQ}){\rightarrow }H_*(X,Y'\cup X_2;{\QQ})\stackrel{(m_*)^{-1}}{\rightarrow }H_*(Y_2,Y'\cup (X_2\cap Y_2);{\QQ})$$
where $m:(Y_2,Y'\cup (X_2\cap Y_2))\into (X,Y'\cup X_2)$ is the inclusion. Clearly, the proposition is a consequence of the following claim. 

\begin{clm}\label{clm diag for cap}
The diagrams
\[\xymatrix{ 
H^*(X';{\QQ})\ar[d]^{(l^*)^{-1}\circ d^*}\ar[rr]^{j_*\circ (\,\,\,\cap \alpha ')}&&
H_*(X,Y;{\QQ})\ar[d]^{(k_*)^{-1}\circ (c_*\circ d_*)}  \\
H^*(X,X_2;{\QQ})\ar[rr]^{\cap \alpha _1\circ b^*}&&
H_*(Y_2,Y';{\QQ})
.}
\] 
and
\[\xymatrix{ 
H^*(X,X_2;{\QQ})\ar[d]\ar[rr]^{\cap \alpha _1\circ b^*}&&
H_*(Y_2,Y';{\QQ})\ar[d]  \\
H^*(X;{\QQ})\ar[rr]^{\cap \alpha }&&
H_*(X,Y';{\QQ})
.}
\] 
are commutative.
\end{clm}

\medskip
\noindent
\textbf{Proof of Claim \ref{clm diag for cap}:} Below, we also denote by $\cap $ the chain map which induces the cap product.We will as well identify cycles (resp., cocycles) with their images under chain maps (resp., cochain maps) induced by certain inclusion maps.

Let $\beta \in H^*(X';{\QQ})$ and take $z\in Z^*(X';{\QQ})$ a representative of $\beta $ and extend it (by zero outside $S_*(X';{\QQ})$) to $z '\in S^*(X;{\QQ})$. Then $\partial ^*z '_{|S_*(X_1;{\QQ})}$ represents $d^*\beta $. By excision axiom, there is $w\in Z^*(X,X_2;{\QQ})$ such that $w_{|S_*(X_1;{\QQ})}=\partial ^*z'_{|S_*(X_1;{\QQ})}+\partial ^*w'$ where $w'\in S^*(X_1,X';{\QQ})$. Extend $w'$ (by zero outside $S_*(X_1;{\QQ})$) to $w''\in S^*(X,X_2;{\QQ})$, and replace $w$ by $w-\partial ^*w''$. The new cocycle $w$ then satisfies $w_{|S_*(X_1;{\QQ})}=\partial ^*z'_{|S_*(X_1;{\QQ})}$ and represents the image of $\beta $ in $H^*(X,X_2;{\QQ})$ and $H^*(X;{\QQ})$.

Because $X_1\cap Y_2$, $X_2\cap Y_1$ and $X'$ are open definable subsets of $X$ which cover $X$, by Lemma \ref{lem omin exci}, we can find a representative $a$ of $\alpha $ such that $a=a_1+a_2+a'$ with $a_1\in S_*(X_1\cap Y_2;{\QQ})$, $a_2\in S_*(X_2\cap Y_1;{\QQ})$, $a'\in S_*(X';{\QQ})$ and, of course, $\partial _*a\in S_*(Y';{\QQ})$. Then $a'$ represents $\alpha '$ and $a_1$ represents $\alpha _1$. It follows that the image of $\beta $ in $H^*(Y,Y_1;{\QQ})$ along the two ways of the first diagram of the claim has representative $w\cap a_1$ and $\partial _*(z\cap a')=(-1)^{{\rm deg}z}z\cap \partial _*a'$ (by Theorem \ref{thm cap product} (6)). We claim that these elements determine the same homology class. Note that, since $z'_{|S_*(X_1\cap Y_2;{\QQ})}=0$ we have $z'\cap a_1=0$. Hence, by Theorem \ref{thm cap product} (6), we have $\partial _*(z'\cap a_1)=\partial ^*z'\cap a_1+(-1)^{{\rm deg}z}z'\cap \partial _*a_1$ $=\partial ^*z'\cap a_1+(-1)^{{\rm deg}z}z'\cap \partial _*a-(-1)^{{\rm deg}z}z'\cap \partial _*a_2-(-1)^{{\rm deg}z}z'\cap \partial _*a'$. But we have (i) $\partial ^*(z'\cap a_1)=w\cap a_1$ since $a_1\in S_*(X_1;{\QQ})$ and $w_{|S_*(X_1;{\QQ})}=\partial ^*z'_{|S_*(X_1;{\QQ})}$; (ii) $z'\cap \partial _*a'=z\cap \partial _*a'$ since $\partial _*a'\in S_*(X';{\QQ})$ and $z'_{|S(X';{\QQ})}=z$ and (iii) $z'\cap \partial _*a-z'\cap \partial _*a_2\in S_*(Y_1;{\QQ})$ because $\partial _*a$ and $a_2$ are in $S_*(Y_1;{\QQ})$. Hence, $\partial _*(z'\cap a_1)-[w\cap a_1-(-1)^{{\rm deg}z}z\cap \partial _*a']\in S_*(Y_1;{\QQ})$ as required.

It remains to show the commutativity of the second diagram of the claim. Let $\gamma \in H^*(X,X_2;{\QQ})$ and take $u\in Z^*(X,X_2;{\QQ})$ a representative of $\gamma $. Then we have $u\cap (a_2+a')=0$ because $a_2+a'\in S_*(X_2;{\QQ})$ and $u_{|S_*(X_2;{\QQ})}=0$. Hence, $u\cap a=u\cap a_1$ as required.
\qed

\qed

\begin{cor}\label{cor extra cap product}
If we take $X_1=X$, $X_2=V$, $Y_1=U$ and $Y_2=W$ in Proposition \ref{prop extra cap product}, then we get the following commutative diagram
\[\xymatrix{ 
H^*(V;{\QQ})\ar[d]^{\cap \alpha '}\ar[rr]^{d^*}&&
H^*(X,V;{\QQ})\ar[rr]^{b^*\circ (l^*)^{-1}}&&
H^*(W,V\cap W;{\QQ})\ar[d]^{k_*\circ (\,\,\,\cap \alpha _1)}  \\
H_*(V,V\cap W;{\QQ})\ar[rr]^{j_*}&&
H_*(X,W;{\QQ})\ar[rr]^{d_*}&&
H_*(W,U;{\QQ})
.}
\] 
\end{cor}

\end{section}

\end{document}